\documentclass[a4paper,12pt]{article}
\usepackage{latexsym,bm}
\usepackage{amsmath,amsfonts,amssymb}
\numberwithin{equation}{section} \numberwithin{prop}{section}
\numberwithin{thm}{section} \numberwithin{lem}{section}
\numberwithin{dfn}{section} \numberwithin{cor}{section}

\title{Matzoh ball soup in spaces of constant curvature
 }

\author{Genqian Liu \\[5pt]
Department of Mathematics\\
 Beijing Institute of Technology\\
Beijing, 100081, P.R. China\null\\[1.5pt]
E-mail: liugqz@bit.edu.cn}

\date{}

\begin{document}
\maketitle
\begin{abstract}

In this paper, we generalize Magnanini-Sakaguchi's result \cite{MS3}
 from Euclidean space to spaces of constant curvature.
 More precisely, we show that if a conductor satisfying the
 exterior geodesic sphere condition in the space of constant
 curvature has initial temperature $0$ and its boundary is kept at
 temperature $1$ (at all times), if the thermal
  conductivity of the conductor is inverse of its metric,
  and if the conductor contains a proper sub-domain,
  satisfying the interior geodesic cone condition and having constant
  boundary temperature at each given time, then the conductor must
  be a geodesic ball.
  Moreover, we show similar result for the
  wave equations and the Schr\"{o}dinger equations
  in spaces of constant curvature.

\end{abstract}
\bigskip\bigskip\bigskip

\section{Introduction}

Klamkin's conjecture \cite{Kl} (also referred to by L. Zalcman in
\cite{Za} as the {\it  Matzoh ball soup problem}) states that, in a
bounded domain
 $\Omega$ (i.e., the {\it Matzoh ball} in ${\Bbb R}^n$), if
 the normalized temperature $u=u(t,x)$ satisfies the heat equation:
  \[ \left\{ \begin{array} {ll}  \frac{\partial u}{\partial t}
  =\triangle u \quad &\mbox{in $\,(0, +\infty)\times \Omega$},\\
  u=1 \quad & \mbox{on $\,(0, +\infty)\times\partial \Omega$},\label {1-1}\\
  u=0 \quad &\mbox{on $\,\{0\}\times\Omega$},\end{array}\right.\]
 and if {\it all} spatial isothermic surfaces of $u$ are {\it invariant with
 time} (the values of $u$ vary with time on its spatial isothermic surfaces),
 then $\Omega$ must be a ball.

 In \cite{A1}-\cite{A2}, this conjecture had been settled affirmatively by
 G. Alessandrini (also see \cite{Sa1} for a different method, by
 which Klamkin's conjecture can be proved).
  A stronger result has also been
 obtained by Magnanini and Sakaguchi in \cite{MS3},
 which says that Klamkin's conjecture holds only
 if {\it one} spatial isothermic surface of $u$ is {\it invariant
 with time}.

  It is a natural question to ask whether Magnanini-Sakaguchi's stronger result
   still  hold in the space ${\Bbb M}_k$
  of constant curvature $k$ ($\, k\in {\Bbb R}^1$)?

\vskip 0.23 true cm

 The main purpose of this paper is to prove the following:

\vskip 0.20 true cm

\noindent  {\bf Theorem 1.1.} \ \  {\it  Let $\Omega$ be a bounded
domain in the $n$-dimensional space ${\Bbb M}_k$ of constant
curvature $k$ with metric $g_{ij}=
 \frac{4\delta_{ij}}{(1+k|x|^2)^2}\,$ ($\delta_{ij}$ is Kronecker's symbol;
 in case of $k>0$, $\Omega$ is
required to lie in a hemisphere), $\,n\ge 2$. Let $\Omega$ satisfy
the exterior geodesic sphere condition and assume that $D$ is a
domain, with boundary $\partial D$, satisfying the interior geodesic
cone condition, and such that $\bar D\subset \Omega$.

Let $u$ be the solution to the problem
  \begin{eqnarray}  \frac{\partial u}{\partial
  t}&=& \sum_{i=1}^n \bigg(\frac{1+k|x|^2}{2}\bigg)^2
  \,\,\frac{\partial^2 u}{\partial x_i^2}
  \quad \;\; \mbox{in $\,(0,+\infty)\times \Omega$}
   \label {1-1},\end{eqnarray}
  and the two conditions:
  \begin{eqnarray}
  u&=&1 \label{1-2}\quad \,  \mbox{on $\,(0,
  +\infty)\times\partial \Omega$},\\
  u&=&0 \label{1-3}\quad \ \mbox{on $\,\{0\}\times \Omega$}.
 \end{eqnarray}
  If $u$ satisfies the extra condition:
   \begin{eqnarray}u(t,x)=a(t), \quad \;
   (t,x)\in (0, +\infty)\times \partial D,\end{eqnarray}
  for some function $a: (0, +\infty)\to (0, +\infty),$
  then $\Omega$ must be a geodesic ball in ${\Bbb M}_k$. }
 \vskip 0.28 true cm

 It is well-known (see, for example, [24, p.79])
  that in a solid medium, the heat flow is governed by two characteristics,
 conductivity and capacity, which may vary over the medium.
 A general mathematical model is provided
  by a manifold $M$, in which the conductivity, or rather
  its inverse, the resistance, corresponds to a Riemannian metric,
  and the capacity corresponds to a Borel measure.
   The above theorem means that in the space of constant
 curvature with the metric $\frac{4\delta_{ij}}{(1+k|x|^2)^2}\,$,
  if the thermal conductivity of the conductor
   is inverse of its metric, and if one spatial
   isothermic surface is invariant with time
   (of course, its boundary is kept at temperature $1$),
  then the conductor takes the shape of a geodesic ball.
 Clearly, Theorem 1.1 reduces to Magnanini-Sakaguchi's result [23] when $k=0$.

 The proof of our main theorem is essentially based on three ingredients:
 The first ingredient is Varadhan's theorem, which not only implies
  that (1.1) is the correct form of the heat equation on ${\Bbb M}_k$,
  but also tells us that $\partial \Omega$ and
 $\partial D$ are equidistant surfaces.
 The second ingredient is a new method which is due to Magnanini
  and Sakaguchi (see \cite{MS3}). This method contains an integral
  transform with respect to time variable, two kinds of balance laws and an asymptotic
  formula.  In order to apply Magnanini-Sakaguchi's method to
  fit our manifold setting, we use two techniques: One is the
  invariance property of operator $\sum_{i,j=1}^n
  (\frac{1+k|x|^2}{2})^2 \frac{\partial^2}{\partial x_i\partial
  x_j}$ under isometries. The other is an orthogonal projection from
  the sphere ${\Bbb S}^n_{1/\sqrt{k}}$ or hyperboloid model ${\Bbb H}^n_{1/\sqrt{-k}}$
  to the Euclidean space
  $\{(x,0)\in {\Bbb R}^{n+1}\big| x\in {\Bbb R}^n\}$, which allows us to
  derives a formula for the principal curvatures (see Lemma 4.1). This is also a key
   step toward the proof of our main theorem.
  The last ingredient is Alexandrov's theorem \cite{A} that provides a
  characteristic property of geodesic spheres in the spaces of
  constant curvature.

   Finally, we show similar result for the
  wave equations and the Schr\"{o}dinger equations
  in spaces of  constant curvature.
 \smallskip

\section{Preliminaries}

\vskip 0.48 true cm

Let ${\Bbb M}_k$ be a complete, simply connected $n$-dimensional
Riemannian manifold of constant curvature $k$.  Then ${\Bbb M}_k$ is
uniquely determined, up to isometric equivalence (see \cite{Hi},
 \cite{Le} or \cite{Wo}). Of course, when $k=0$ we may take ${\Bbb M}_k={\Bbb
R}^n$ with the usual Euclidean metric $ds^2 = dx_1^2 +dx_2^2 +\cdots
+dx_n^2$. For $k\ne 0$, various realizations
 are  possible. Thus, if $k>0$ we may take for ${\Bbb M}_k$ the
 $n$-dimensional sphere ${\Bbb S}^n_\rho=\{y\in {\Bbb R}^{n+1}\big| \sqrt{y_1^2 +\cdots
 +y_{n+1}^2}=\rho\}$ of radius $\rho=1/\sqrt{k}$,
 centered at the origin in ${\Bbb R}^{n+1}$, with the induced
 Euclidean metric; equivalently, ${\Bbb S}^n_\rho$ may be realized
  by stereographic projection from the north pole.
  This is a map
  $\sigma: {\Bbb S}^n_\rho\setminus \{(0,\cdots, 0,\rho)\}\ni y\mapsto x=\sigma y\in {\Bbb
  R}^n$, which maps a point $y\in {\Bbb S}^n_\rho$ into the intersection $x\in
  {\Bbb R}^n$ of the line jointing $y$ and the north pole $(0, \cdots,
  0,\rho)$ with the equatorial hyperplane ${\Bbb R}^n$.
  Clearly, the south pole $(0,\cdots, 0, -\rho)$ is mapped into the
  origin, and one has (see [8, p.59])
  $$ y_{n+1} =\rho\,\frac{|x|^2-\rho^2}{|x|^2+\rho^2}, \quad\; \,  (y_1, \cdots,
  y_n)=\frac{2\rho^2 x}{\rho^2+|x|^2}\quad\;
  \left(x=\frac{\rho}{\rho-y_{n+1}}(y_1, \cdots, y_n)\right).$$
    The map $\sigma$ induces a matric on ${\Bbb R}^n$:
 \begin{eqnarray} ds^2 =\frac{4|dx|^2}{(1+ |x|^2/\rho^2)^2},\end{eqnarray}
i.e.,
  \begin{eqnarray} g_{ij}=\bigg\langle \frac{\partial}{\partial x_i},
  \frac{\partial}{\partial x_k}\bigg\rangle =\frac{4\delta_{ij}}
  {\big(1+|x|^2/ \rho^2\big)^2}.\end{eqnarray}

  Let $M$ and $N$ be two manifolds with metrics $g$ and $h$, respectively.
    We say that a diffeomorphism $\Phi: (M, g)\to (N, h)$ is an
  isometry if $\Phi^* h=g$.
  It is well-known that every isometry of ${\Bbb S}^n_\rho$ is an
  element of $O(n+1)$.  For $y\in {\Bbb S}^n_\rho\setminus \{\mbox{the north
  pole}\}$, take $R_y\in O(n+1)$ satisfying $R_y(y)=\mbox{the south
  pole}\,$  ($R_y (-y)=\mbox{the north pole}$).
  Setting $z=\sigma y\in{\Bbb R}^n$ and $z^*:=\rho^2
 z/|z|^2$ yields $\sigma (-y)=-z^*$. Let us consider the map $f=\sigma\circ
R_y \circ \sigma^{-1}$.
   One has $f(z)=0$, $f(-z^*)=\infty$.
 Note that a M\"{o}bius transformation with such property is
 of the form
 $$f(x)=\lambda A((x+z^*)^* -(z+z^*)^*)$$ with $\lambda>0$ and a
 constant orthogonal matrix $A$ (see [Ah, p.21]). Similar to
 the method of [5, p.1106], we get that
 \begin{eqnarray} f(x)= \frac{\rho^2(\rho^2+|z|^2)}{|z|^2}
 A((x+z^*)^* - (z+z^*)^*) \mbox{
with }A\in O(n).\end{eqnarray}

\vskip 0.26 true cm

  Let ${\Bbb R}^{n+1}$ be equipped with the Lorentzian metric
    $$\langle y,y\rangle = -y_{n+1}^2 +y_1^2+\cdots y_n^2.$$
For $k<0$, let $\rho=1/\sqrt{-k}$ and $${\Bbb H}^n_\rho =\{y\in
{\Bbb R}^{n+1}\big|\langle y,
 y \rangle=-\rho,\;  y_{n+1}>0\}$$
  with the Riemannian metric induced from the Lorentzian metric.
 ${\Bbb H}_\rho^n$ is called the hyperboloid model or
  Lobochevskian pseudo-sphere (see [18, p.38-42]).
  By regarding ${\Bbb R}^n$ as $\{(x,0)\in {\Bbb
  R}^{n+1}\}$,
   we consider the {\it hyperbolic stereographic projection}
   $\zeta: {\Bbb H}^n_\rho\ni y\to x=\zeta y\in {\Bbb R}^n$, which map a
   point $y\in {\Bbb H}^n_\rho$ into the intersection $x\in {\Bbb
   B}_\rho:= \{x\in {\Bbb R}^n \big| |x|< \rho\}$ of line joining
   $y$  and the point $(0, \cdots, -\rho)$ with ${\Bbb R}^n$.
   Then, the point $(0, \cdots, 0, \rho)$ is mapped into the
    origin, and we have
    $$ y_{n+1}= \rho\,\frac{\rho^2+|x|^2}{\rho^2-|x|^2},
    \quad  \; (y_1, \cdots, y_n)= \frac{2\rho^2 x}{\rho^2-|x|^2} \;
    \; \quad\;  \; \big(x=\frac{\rho}{\rho +y_{n+1}}(y_1, \cdots, y_n)\big).$$
This map induces the metric on ${\Bbb B}_\rho$:
   \begin{eqnarray} ds^2 =\frac{4 |dx|^2}{(1-|x|^2/\rho^2)^2},\end{eqnarray}
i.e.,
  \begin{eqnarray} g_{ij}=\bigg\langle \frac{\partial}{\partial x_i},
  \frac{\partial}{\partial x_k}\bigg\rangle =\frac{4\delta_{ij}}
  {(1-|x|^2/ \rho^2)^2}.\end{eqnarray}
  Conformal transformations of the plane are holomorphic mappings,
 whereas  in higher dimensions ($n\ge 3$) the only possibilities
 are rotations, dilations, inversions $x\to x^{*}=\frac{\rho^2 x}{|x|^2}$
 and their compositions (Liouville's theorem, see [10, \S 15]).
  Every isometry of $({\Bbb B}^n_\rho, g)$ is a conformal map $f:{\Bbb B}^n_\rho\to
 {\Bbb B}^n_\rho$.
   Similar to \cite{Ah}, we can verify that the general form of such a map is:
  \begin{eqnarray} f(x)=T_z(x):=
  A\frac{\rho^2[(\rho^2-|z|^2)(x-z)-|x-z|^2z]}
  {\rho^4+|x|^2|z|^2 -2\rho^2 x\cdot z} \;\; \mbox{
with }A\in O(n).\end{eqnarray}
 Obviously, $T_z(z)=0$, and the isometries of $({\Bbb B}^n_\rho, g)$
 transform spheres into spheres.

 Throughout this paper (except for the proofs of Lemma 4.1 and Theorem 4.2),
  ${\Bbb M}_k$ can be regarded as ${\Bbb R}^n$
  with metric (2.2) when $k>0$;
  ${\Bbb M}_k={\Bbb R}^n$
  with the Euclidean metric when $k=0$; ${\Bbb M}_k$ as
 ${\Bbb B}^n_\rho$ with metric (2.5) when $k<0$.
  As $k\to 0$, the metric (2.1) and (2.4) approach the flat
   (Euclidean) metric. For (2.1) this is geometrically obvious,
   and in any case can be seen from the equivalent form
 \begin{eqnarray} ds^2_k
   =\frac{4|dx|^2}{(1+k|x|^2)^2},\end{eqnarray}
valid for all $k$.
 The pair $(D, ds_k^2)$ with $D\subset {\Bbb R}^n$ is call the
 {\it canonical representation} of the domain $D$ in ${\Bbb M}_k$.

 \vskip 0.28 true cm

Let $M$ be an $n$-dimensional Riemannian manifold with the metric
 $g_{ij}(x)$, and let $L$ be the following differential operator acting on smooth
 functions on $M$:
   $$Lu=\frac{1}{2}\sum_{i,j=1}^n g^{ij} (x) \frac{\partial^2 u}{\partial x_i
   \partial x_j},$$
 where $(g^{ij} (x))$ is the matrix inverse to $(g_{ij}(x))$.

 Let $p(\tau)$, $0\le \tau \le 1$, be a  smooth path in $M$.
 Then the length of such a path is defined as
  $$ l(p) = \int_0^1 [\dot{p}(\tau) g(p(\tau)) \dot{p}(\tau)]^{1/2}
  d\tau,$$
  where $\dot{p}(\tau)$ stands for $dp(\tau)/d\tau$
  and $(\theta g\theta)$ for the quadratic
  form $\sum_{i,j=1}^n  g_{ij} (x)\theta_i \theta_j$;
  $\,l(p)$ is the natural length in a metric defined locally as
  $$ ds^2 =\sum_{i,j=1}^n  g_{ij} dx_i dx_j.$$
The  global distance $d(x, y)$ induced by this metric is defined as
\begin{eqnarray*} d(x,y)= \inf_{\{p\big|\,p(0)=x, \, p(1)=y\}}\,
l(p).\end{eqnarray*}

\vskip 0.23 true cm

 \noindent  {\bf Lemma 2.1 (Varadhan's theorem, see \cite{Va}).} \ \ {\it Let $\Omega$
  be a bounded domain in Riemannian manifold $M$
  with uniform H\"{o}lder continuous metric $g_{ij}(x)$.
       Let $\phi(s,x)$ be the solution of the equation
     \begin{eqnarray} \frac{1}{2} \sum_{i,j=1}^n  g^{ij}(x)
     \frac{\partial^2 \phi}{\partial x_i
     \partial x_j}=s\phi
     \quad \mbox{for}\;\; x\in \Omega\end{eqnarray}
with the boundary value $\phi=1$ on the boundary $\partial \Omega$
of $\Omega$.
 Then
 \begin{eqnarray} \lim_{s\to \infty} \left[-\frac{1}{\sqrt{2s}} \log \,\phi(s, x)\right] =
 \mathfrak{F}(x), \end{eqnarray} uniformly over compact subset
 of $\bar \Omega$,
where $x$ is any point of $\bar \Omega$ and \begin{eqnarray}
\mathfrak{F}(x)= \mbox{dist}\,(x, \partial \Omega)\end{eqnarray}
 is the shortest distance to the boundary $\partial \Omega$ from $x$}.

\vskip 0.26 true cm

 \noindent {\bf Lemma 2.2 (Alexandrov's theorem,
see \cite{A}).} \ \  {\it Let $\Gamma$ be a closed
$(n-1)$-dimensional surface in an $n$-space ${\Bbb M}_k$ of constant
curvature $k$ (in case of sphere, $\Gamma$ is required to lie in a
hemisphere). Suppose that $\Gamma$ has no multiple points and is of
class $C^2$.
 Let $\lambda_1\ge \cdots \ge \lambda_{n-1}$ denote its principal
curvatures, at an arbitrary point $p\in \Gamma$. Assume that
$F=F(\beta_1, \cdots, \beta_{n-1})$ is a continuous differentiable
function, defined for $\beta_1, \cdots,
 \beta_{n-1}$, and subject to the condition
 \begin{eqnarray*} \text{const}> \frac{\partial
  F(\beta_1, \cdots, \beta_{n-1})}{\partial \beta_j}>\text{const}> 0,
  \quad \,(j=1, \cdots, n-1),\end{eqnarray*}
 at least on $\Gamma$, i.e., $\beta_j=\lambda_j$  $\; (j=1, \cdots,
 n-1)$.  Then, if $F(\lambda_1, \cdots, \lambda_{n-1})\equiv \text{constant on }$
 $\Gamma$, $\,\,\Gamma$ is a geodesic sphere}.

\vskip 0.25 true cm

 \noindent {Proof.}  \ When $n=2$, we have that
 $$const>\frac{dF(\beta_1)}{d\beta_1}>const>0\quad \; \mbox{on}\;\; \Gamma,$$ i.e.,
   $F(\lambda_1)$ is increasing in $\lambda_1\in\Gamma$.
   Thus, from $F(\lambda_1)\equiv constant$ on $\Gamma$, we get
   that $\lambda_1$ (i.e., the curvature of $\Gamma$)
    must be a constant on $\Gamma$, which implies that
   $\Gamma$ is the boundary of a geodesic disk in ${\Bbb M}_k$.
 When $n \ge 3$, the theorem had been  proved by A. D. Alexandrov
 (see [2, Theorem and ($I_2$) of Remark (6)]). $\;\;\square$

 \vskip 1.32 true cm

\section{Isometric invariance and balance law}

\vskip 0.46 true cm

In this section, we shall prove some lemmas, which are needed for
proving our main theorem.
 First, we prove a simple invariance property of the operator (3.1) below.
  If $(U, \phi)$ is a local chart on $M$ and $f\in C^2(M)$, we often
  write $f^*$ for the composite function $f\circ \phi^{-1}$.

\vskip 0.2 true cm

 \noindent{Lemma 3.1.} \ \  {\it Let $\Phi$ be a diffeomorphism
  of the Riemannian manifold $M$ with metric $g_{ij}$.
 Then $\Phi$ leaves the operator
 $L$ invariant if and only if
 it is an isometry,
  where \begin{equation}\label{3-1}  L f=\sum_{i,j=1}^n g^{ij} (x)
  \frac{\partial^2 f}{\partial x_i\partial
  x_j}.\end{equation}}
 \vskip 0.25 true cm

 \noindent{Proof.}   Let $p\in M$ and let $(V,
\psi)$ be a local chart
 around $p$. Then $(\Phi(V), \psi\circ \Phi^{-1})$ is a local chart
 around $\Phi(p)$. For $x\in V$, let $y= \Phi(x)$ and
   \begin{eqnarray*} &\quad \;\; \psi(x)=(x_1, \cdots, x_n), \quad \; x\in V,\\
   &(\psi\circ \Phi^{-1})(y)=(y_1, \cdots, y_n), \quad \; y\in
   \Phi(V).\end{eqnarray*}
  Then $$x_i(x)=y_i(\Phi(x)), \quad \, d\Phi_x
  \left(\frac{\partial}{\partial x_i}\right)_x
  =\left(\frac{\partial}{\partial y_i}\right)_{\Phi(x)} \quad \,
  (1\le i\le n),$$
  where $d\Phi_x$ is the tangent map.
   For each function $f\in C^2(M)$,
   \begin{eqnarray} ((Lf)^{\Phi^{-1}})(x)= (Lf)(\Phi(x))=
     \sum_{i,j=1}^n g^{ij}(y) \frac{\partial^2 f^*}{\partial y_i
     \partial y_j},\end{eqnarray}
\begin{eqnarray} (Lf^{\Phi^{-1}})(x)=
     \sum_{i,j=1}^n g^{ij}(x) \frac{\partial^2 (f\circ \Phi)^*}{\partial x_i
     \partial x_j}.\end{eqnarray}
 Now if $\Phi$ is an isometry, then $g_{ij}(x)=g_{ij}(y)$ for all $i,
  j$. Because of the choice of coordinates,  we have
 $$\frac{\partial^2 f^*}{\partial y_i \partial y_j}= \frac{\partial^2 (f\circ
 \Phi)^*}{\partial x_i \partial x_j},\quad \; \, 1\le i, j\le n.$$
    Thus the right-hand sides of (3.2) and (3.3) coincide and
  $L^\Phi=L$, which implies that the operator
  $L$ is invariant.
   On the other hand, if (3.2) and (3.3) agree, then we find by
   equating coefficients that $g_{ij}(x)=g_{ij}(y)$, which shows that
   $\Phi$ is an isometry.  $\;\; \square$

\vskip 0.26 true cm

 The following Lemma is the so-called {\it balance law}, which has
 been proved by Magnanini and Sakaguchi in Euclidean case (see
 \cite{MS1}, \cite{MS2}, \cite{MS3}) and by Sakaguchi
 in ${\Bbb M}_k$ with
  the Laplace-Beltrami operator instead of $L$ (see \cite{Sa2}).

\vskip 0.26 true cm

  \noindent{Lemma 3.2.} \ \  {\it Let $\Omega$ be a domain
  in the $n$-dimensional space  ${\Bbb M}_k$ of constant
  curvature $k$ (in case of sphere, $\Omega$ is required
  to lie in a hemisphere), $n\ge 2$.
    Let $x_0$ be a point in
  $\Omega$ and set $d_* =\text{dist}\, (x_0, \partial \Omega)$.
  Assume that $v=v(t, x)$ is a solution of
  \begin{eqnarray}\frac{\partial v}{\partial t}=\sum_{i,j=1}^n
  \bigg(\frac{1+k|x|^2}{2}\bigg)^2 \,\frac{\partial^2 v}{\partial x_i\partial x_j} \,
  \quad \mbox{in} \;\; (0,+\infty) \times \Omega.\end{eqnarray}
Then, the following two assertions hold:

  (i) \ \ $v(t, x_0)=0$  for every $t\in (0, +\infty)$ if and only if
\begin{equation}  \int_{\partial B_r(x_0)}   v(t,x)dA_r
=0\label{3-1}\quad \;  \mbox{for every $(t, r)\in (0,+\infty)\times
[0, d_*)$},\end{equation}
     where $\partial
     B_r(x_0)$ denotes the geodesic sphere centered at $x_0$ with
     radius $r>0$ and $dA_r$ denotes its area element;

  (ii)  \ \  $\nabla \, v(t, x_0)=0$  for every $t\in (0, +\infty)$ if and only if
\begin{equation}  \int_{\partial B_r(x_0)}  \exp_{x_0}^{-1} x v(t,x)dA_r
=0\label{3-2}\quad \;  \mbox{for every $(t, r)\in (0,+\infty)\times
[0, d_*)$},\end{equation}}
 where $\exp_{x_0}$ is the exponential map at $x_0$.

\vskip 0.28 true cm

\noindent {\it Proof.} \ \  (i) \ \  If (3.5) holds, then we
immediately get that $v(t, x_0)=0$ for
 every $t\in (0, +\infty)$. Conversely,
  for any two points $x'$ and $x''$ in ${\Bbb M}_k$, we can find an
  isometry $\Phi$ that maps ${\Bbb M}_k$ onto itself such that $\Phi x'=x''$
  (cf. section 2).
  It follows from Lemma 3.1 that the operator $L$ is invariant under
  the isometry $\Phi$ (here $g^{ij}(x)= \big(\frac{1+k|x|^2}{2}\big)^2 \delta_{ij}$),
  that is, $(Lu)\circ \Phi=L(u\circ \Phi)$.
 Thus, by an isometry we may put $x_0=0$ in the canonical
 representation.
  Note that spherical coordinates are valid about {\it any point} in
  $\Omega\subset {\Bbb M}_k$ for each fixed $k$ (see [7,
  p.37-39]). Therefore, about the origin in the canonical representation,
  there exists a coordinate system $(r, \theta)\in [0, d_*)\times {\Bbb
 S}^{n-1}$, relative to which the Riemannian metric reads
 as
 \begin{eqnarray}  ds^2 = (dr)^2 + (h_k(r))^2  |d\theta|^2,\end{eqnarray}
 where \begin{eqnarray}  h_k(r)=\left\{ \begin {array}{ll} (1/\sqrt{k})\sin \sqrt{k} \,r,
  & \quad k>0,\\
  r, & \quad k=0,\\
  (1/\sqrt{-k})\sinh {\sqrt{-k}}\,r, & \quad k<0,\end{array}\right.\end{eqnarray}
  $|d\theta|^2$ denotes the metric on the Euclidean sphere ${\Bbb S}^{n-1}$
  of radius $1$, and $r$ the geodesic distance from $x_0=0$.

   Let $C(r):=\{\theta\in T_0 {\Bbb M}_k\big||\theta|=1\;\; \mbox{and}\;\;
   \gamma_\theta (s)=\exp_0(s\theta),\, s\in [0,r]$, is minimizing$\,\}$.
  In view of $T_0{\Bbb M}_k ={\Bbb R}^{n}$, we see that  $C(r)= {\Bbb S}^{n-1}$
   for all $r\in[0, d_*)$.
  Denote $x=x(r, \theta)\in {\Bbb M}_k$, where $dist \,(x,0)=r$, $\,\theta\in
{\Bbb S}^{n-1}$.
   It follows that
   \begin{eqnarray*} \int_{B_r(0)} v(t, x) d\mu(x)=
   \int_0^r \left(\int_{{\Bbb S}^{n-1}} v(t, \exp_0(\bar r, \theta))
   J(\bar r, \theta)d\Theta(\theta)\right)d\bar r, \end{eqnarray*}
which implies
 \begin{eqnarray} \int_{\partial B_r(0)} v(t, x) dA_r&=&
  \int_{{\Bbb S}^{n-1}} v(t, \exp_0(r,\theta))
   J(r, \theta)d\Theta(\theta)\\
  &=&\int_{{\Bbb S}^{n-1}} v(t, \exp_0(r \theta))\,
   (h_k(r))^{n-1} d\Theta(\theta),\nonumber \end{eqnarray}
where $d\Theta$ is the volume form of the unit $(n-1)$-sphere, $J(r,
\theta)=\sqrt{\mbox{det}(g_{ij})}=(h_k(r))^{n-1}$ (see [9,
p.74-76]). Then (3.5) is equivalent to
\begin{eqnarray*} \int_{{\Bbb S}^{n-1}} v(t, \exp_0(r\theta)) d\Theta(\theta)=0 \quad
\mbox{for any } (t, r)\in (0, +\infty)\times [0,d_*).\end{eqnarray*}
We define  \begin{eqnarray}  U(t, r):=\int_{{\Bbb S}^{n-1}} v(t,
\exp_0 (r\theta))d\Theta(\theta), \end{eqnarray} for all $(t,r)\in
(0, +\infty)\times [0, d_*)$.
 Since  \begin{eqnarray} L&=&\sum_{i=1}^n \left(\frac{1+k|x|^2}{2}\right)^2
 \, \frac{\partial^2 }{\partial x_i^2}\\
  &=& \left(\frac{1+k r^2}{2}\right)^2\left(\frac{\partial^2}{\partial
  r^2}+\frac{n-1}{r} \, \frac{\partial}{\partial r} +\frac{1}{r^2}
  \triangle_{S^{n-1}}\right)\nonumber, \end{eqnarray}
  where $\triangle_{S^{n-1}}$ denotes  the Laplace-Beltrami operator
  on $S^{n-1}$, by substituting this into
  $$ 0=\int_{{\Bbb S}^{n-1}}
  \left[\big(\frac{\partial}{\partial t}-L\big) v(t, \exp_0(r\theta))\right]d\Theta(\theta),$$
and using $\int_{{\Bbb S}^{n-1}} \big(\triangle_{{\Bbb S}^{n-1}}
v(t,
 \exp_0 (r\theta))\big) d\Theta(\theta) =0$,
  we obtain
  \begin{eqnarray} \frac{\partial U}{\partial t} =\left(\frac{1+kr^2}{2}\right)^2
 \left(\frac{\partial^2 U}{\partial
  r^2}+\frac{n-1}{r} \, \frac{\partial U}{\partial r}\right)\quad \;
  \mbox{in $[0, d_*)\times (0, +\infty)$}. \end{eqnarray}
    It follows from the local regularity result of parabolic
    equations (see [12, 13], [19, 20] or [26, p.404-405])
    that $U$ and $\frac{\partial U}{\partial r}$
    are real analytic in $(0, +\infty)\times [0, d_*)$ .
Therefore
\begin{eqnarray} 4 r \frac{\partial U}{\partial t}=
 \big(k^2 r^5 + 2k r^3 +r\big) \frac{\partial^2 U}{\partial r^2}
 +(n-1) (k^2r^4+ 2kr^2 +1)\frac{\partial U}{\partial r}\end{eqnarray}
 for all  $(t, r)\in (0, +\infty)\times [0, d_*)$.
 Obviously, $U(t, 0)=0$ for any $t>0$. From (3.12), we have
  $$\frac{\partial U}{\partial r} (t,0)=\lim_{r\to 0} \; \frac{r}{n-1}
  \left[\bigg(\frac{2}{1+kr^2}\bigg)^2
   \,\frac{\partial U}{\partial t} -\frac{\partial^2 U}{\partial
   r^2}\right]=0$$
for all $t\in (0, +\infty)$.

We shall show by induction that
   \begin{eqnarray}\frac{\partial^m U}{\partial r^m}(t,0) =0 \quad \; \mbox{for
   all $t>0$ and any integer $m\ge 0$}.\end{eqnarray}
   Suppose that
   $$U(t, 0)=\frac{\partial U}{\partial r}(t, 0)=\cdots
   =\frac{\partial^m U}{\partial r^m}(t, 0)=0 \quad \mbox{for all }
   t >0.$$
By differentiating both sides of (3.13) for $m$ times with respect
to $r$ , we get
   \begin{eqnarray*}  &&\sum_{j=0}^m C_m^j  \;\frac{\partial^j(4r)}{\partial
   r^j}\;
\left( \frac{\partial}{\partial t} \frac{\partial^{m-j} U}{\partial r}\right)\\
  &&\qquad \; =
 \sum_{j=0}^m C_m^j \frac{\partial^j (k^2r^5 + 2kr^3 +r)}{\partial r^j}
 \; \frac{\partial^{m+2-j} U}{\partial r^{m+2-j}}\\
 && \qquad \; \; \,+(n-1)
\sum_{j=0}^m C_m^j \frac{\partial^j (k^2r^4+ 2kr^2 +1)}{\partial
r^j}
 \; \frac{\partial^{m+1-j} U}{\partial r^{m+1-j}},\end{eqnarray*}
where $C_m^j=\frac{m!}{j! (m-j)!}$.
  Thus, letting $r=0$ and using the above assumption, we have
  \begin{eqnarray*} 0= m\, \frac{\partial^{m+1} U}{\partial r^{m+1}}
  + (n-1)\, \frac{\partial^{m+1} U}{\partial r^{m+1}},\end{eqnarray*}
i.e, $\frac{\partial^{m+1} U}{\partial r^{n+1}}(t, 0)=0$. It follows
from induction that (3.14) holds. From the analyticity of $U$,  we
 obtain  that
$$ U\equiv 0 \quad \; \mbox{in} \;\; (0,+\infty)\times [0, d_*).$$
 Therefore, we conclude that (3.5) is true.

 \vskip 0.28 true cm

 (ii)  \ \  As in the argument of (i), by putting $x_0=0$ we get that (3.6) is
 equivalent to
 \begin{eqnarray*} \int_{{\Bbb S}^{n-1}} (r\theta)\,
 v(t, \exp_0 (r\theta))(h_k(r))^{n-1} d\Theta(\theta)=0
 \quad \; \mbox{for all} \;\; (t, r)\in (0, +\infty)\times [0,
 d_*),\end{eqnarray*}
i.e.,
 \begin{eqnarray*} \int_{{\Bbb S}^{n-1}} \theta \,v(t, \exp_0 (r\theta)) d\Theta(\theta)=0
 \quad \; \mbox{for all} \;\; (t, r)\in (0, +\infty)\times [0,
 d_*).\end{eqnarray*}
 If (3.6) holds, then, by the divergence theorem, we get that  $\nabla v(t,0)=0$ for every $t>0$.
 We shall prove the converse assertion.
  Let us introduce an ${\Bbb R}^n$-valued function $Q(t, r)$ by
  \begin{eqnarray} Q(t, r)=\int_{{\Bbb S}^{n-1}} \theta v(t,
  \exp_0(r\theta))d\Theta(\theta) \quad \;\; (t,r)\in (0, +\infty)\times [0,
  d_*).\end{eqnarray}
 By putting (3.11) into
 \begin{eqnarray*} 0= \int_{{\Bbb S}^{n-1}} \theta\left[\left(\frac{\partial}{\partial t}-
  L\right) v(t, \exp_0(r\theta))\right] d\Theta(\theta), \end{eqnarray*}
 and using $-\triangle_{{\Bbb S}^{n-1}} \theta =(n-1) \theta$ together with
 integration by parts, we obtain that in $(0, +\infty)\times [0,
  d_*)$,
\begin{eqnarray} \frac{\partial Q}{\partial t} =\left(\frac{1+kr^2}{2}\right)^2
 \left(\frac{\partial^2 Q}{\partial
  r^2}+\frac{n-1}{r} \, \frac{\partial Q}{\partial r}
  - \frac{n-1}{r^2} Q\right). \end{eqnarray}
   Thus
\begin{eqnarray}   4 r^2 \frac{\partial Q}{\partial t}&=&
 \big( k^2r^6 + 2kr^4 +r^2\big) \frac{\partial^2 Q}{\partial r^2}\\
   && +(n-1) (k^2r^5+ 2kr^3 +r)\frac{\partial Q}{\partial r}
     -(n-1)Q\nonumber \end{eqnarray}
 for all  $(t, r)\in (0, +\infty)\times [0, d_*)$.
  In view of $\nabla v(t,0)=0$,
  we find by the divergence theorem that $Q(t,0)=\frac{\partial Q}{\partial r}(t,0)=0$.
   It follows from the method of
   induction that $\frac{\partial^m Q(t, 0)}{\partial r^m}=0$ for all
 $t\in (0, +\infty)$ and $m=1,2, \cdots$. Therefore,  the analyticity of $Q(t, r)$ implies
  that $Q(t, r)\equiv 0$ in $(0, +\infty)\times [0, d_*)$, and
 the desired result holds.    $\,\, \square$

\vskip 0.28 true cm

  \noindent{Lemma 3.3.} \ \  {\it Let $\Omega$
  be a domain with $C^2$ boundary in the $n$-dimensional space
  ${\Bbb M}_k$ of constant curvature $k$, $n\ge 2$, and let
  $W(s, x)$ be the solution of the following elliptic boundary value problem
    \begin{eqnarray}  \sum_{i=1}^n  \bigg(\frac{1+k|x|^2}{2}\bigg)^2\;
    \frac{\partial^2 W}{\partial x_i^2}=sW \quad \; &\mbox{in $\;\;\Omega$}, \\
     W=1  \qquad \qquad \qquad  \qquad\,\qquad \,
     &\mbox{on $\,\partial \Omega$}.\end{eqnarray}
  Then, for every $\epsilon >0$, there exists a positive number
  $s_\epsilon$ such that
  \begin{eqnarray}  W_\epsilon^{-}(s, x)\le W(s, x)\le
  W_\epsilon^+(s, x)\end{eqnarray}
   for every $x\in \bar \Omega$ and every $s\ge s_\epsilon$, where
   \begin{eqnarray} W_\epsilon^\pm (s, x)=\exp
   \{-\sqrt{s(1\mp\epsilon)}\, \mathfrak{F}(x)\},\end{eqnarray}
  and $\mathfrak{F}(x)$ is defined by (2.10)}.

\vskip 0.23 true cm

 \noindent {\it Proof.} \ \  We can take $\delta>0$  small enough such that
 the function $\mathfrak{F}=\mathfrak{F}(x)$
 defined in (2.10) is of class $C^2$ in the set
 $\overline{\Omega_\delta}$ where
 \begin{eqnarray} \Omega_\delta =\{ x\in \Omega\, : \,
 \mathfrak{F}(x)<\delta\}. \end{eqnarray}

   It is easy to calculate
 \begin{eqnarray*} && \left(\sum_{i=1}^n  \bigg(\frac{1+k|x|^2}{2}\bigg)^2\,
    \frac{\partial^2 W_\epsilon^\pm}{\partial x_i^2}\right) -s W_\epsilon^\pm \\
   &=& \big(\exp\{-\sqrt{s(1\mp\epsilon)}\,\mathfrak{F}(x)\}\big) \left\{
    \bigg(\frac{1+k|x|^2}{2}\bigg)^2 \right.\\
     && \left.\times \sum_{i=1}^n \left[
   -\sqrt{s(1\mp\epsilon)}
   \frac{\partial^2 \mathfrak{F}}{\partial x_i^2} +s(1\mp\epsilon)
   \bigg(\frac{\partial \mathfrak{F}}{\partial x_i}\bigg)^2\right] -s\right\}\\
    &=&\mp\epsilon  \sqrt{s}\,
  \left\{\sqrt{s} \pm\frac{\sqrt{(1\mp\epsilon)}}{\epsilon}\sum_{i=1}^n
  \bigg(\frac{1+k|x|^2}{2}\bigg)^2 \frac{\partial^2 \mathfrak{F}}{\partial x_i^2}\right\}
  W_\epsilon^\pm \quad \quad \; \mbox{in $\Omega_\delta$}. \end{eqnarray*}
   Here we have used the fact that
  \begin{eqnarray*}    \sum_{i=1}^n  \bigg(\frac{1+k|x|^2}{2}\bigg)^2 \bigg(\frac{\partial
   \mathfrak{F}}{\partial x_i}\bigg)^2 = \sum_{i=1}^n
   \frac{\partial \mathfrak{F}}{\partial x_i} \left[\bigg(\frac{1+k|x|^2}{2}\bigg)^2
    \delta_{ij}\right]\frac{\partial
   \mathfrak{F}}{\partial x_j} =\langle \nabla \mathfrak{F},
   \nabla \mathfrak{F}\rangle =1.\end{eqnarray*}
   Setting $M_\delta =\max_{\bar \Omega_\delta} \big|\sum_{i=1}^n
   \left(\frac{1+k|x|^2}{2}\right)^2 \frac{\partial^2 \mathfrak{F}}{\partial x_i^2}\big|$,
   we get that if $s\ge \frac{1+\epsilon}{\epsilon^2} M_\delta^2$,
   then in $\Omega_\delta$
   \begin{eqnarray} &&  \sum_{i=1}^n \bigg(\frac{1+k|x|^2}{2}\bigg)^2
   \, \frac{\partial W_\epsilon^+}
   {\partial x_i^2}  - s
   W_\epsilon^+\le 0\\  &&
   \sum_{i=1}^n \bigg(\frac{1+k|x|^2}{2}\bigg)^2 \, \frac{\partial
   W_\epsilon^-} {\partial x_i^2} -s
   W_\epsilon^{-}\ge 0. \end{eqnarray}
  Since the function $-\frac{1}{\sqrt{s}} \log W(s, x)$ converges
  uniformly on $\bar \Omega$ to $\mathfrak{F}(x)$ as $s\to +\infty$, by Lemma
  2.1 (Varadhan's theorem) there exists a real number $s^*>0$ such that
 for every $s\ge s^*$,
 $$ -\delta (1-\sqrt{1-\epsilon}) \le -\frac{1}{\sqrt{\epsilon}} \log
 W(s, x)-\mathfrak{F}(x) \le \delta (\sqrt{1+\epsilon} -1),   \quad x\in
  \bar \Omega.$$
  Put $s_\epsilon=\max \{s^*,
 \frac{1+\epsilon}{\epsilon^2}M_\delta^2 \}.$
 Completely similar to [23, p.938], we get (3.23). $\;\,\square$

 \vskip 1.32 true cm

\section{Principal curvatures and asymptotic formulas}

\vskip 0.43 true cm

We introduce some notations and definitions for the principal
curvatures of $\partial \Omega$.
 Let $M$ be an $m$-dimensional submanifold  of the $n$-dimensional
 Riemannian manifold $N$. The metric $\langle \cdot, \cdot\rangle$
  on $N$ induces a metric on $M$.  Then one has
  $$\nabla_X^M Y=(\nabla_X^N Y)^\top \quad \, \mbox{for} \;\; X, Y\in
\Gamma(TM),$$ where $\nabla^N$ is the Levi-Civita connection of $N$,
$\nabla^M$ is the induced connection, and $\top: T_x N\to T_x M$ for
$x\in M$ denotes the orthogonal projection.

 Let $\nu(x)$ be a vector field in a neighborhood of $x_0\in M\subset
 N$, that is orthogonal to $M$,
  i.e.,
  \begin{eqnarray} \langle \nu(x), X\rangle =0 \quad \, \mbox{for all $X\in T_x
  M$}.\end{eqnarray}
  We denote by $T_xM^\bot$ the orthogonal complement of $T_x M$ in $T_x
  N$. The bundle $TM^\bot$ with fiber $T_x M^\bot$ at $x\in M$ is called normal bundle of
  $M$ in $N$. (4.1) means $\nu(x)\in T_x M^\bot$.
 For a fixed normal field $\nu(x)\in T_x M^\bot$, we write $A_\nu
 (X)=(\nabla^N_X \nu)^\top$. Clearly, $A_\nu : T_x M \to T_x M$
 is selfadjoint with respect to the metric $\langle \cdot,
 \cdot\rangle$.  Suppose $\langle \nu(x), \nu(x)\rangle \equiv 1$;
 i.e., $\nu$ is a unit normal field.   The $m$ eigenvalues  of
 $A_\nu$ which
 are all real by self adjointness are called the {\it principal
 curvatures} of $M$ in the direction $\nu$, and the corresponding
 eigenvectors are called principal curvature vectors.

   For any  point $x\in \bar \Omega\subset {\Bbb M}_k$,
   let $\mathfrak{F}(x)$ be defined by (2.10). Then $\mathfrak{F}(x)=0$ is the
   hypersurface $\partial \Omega$.
        Since $\nu(x)=\nabla \mathfrak{F}(x)$ for any $x\in \partial \Omega$, we know that
  $\nabla_X^{\Omega} \nabla \mathfrak{F}(x)$ is always tangential to $\partial \Omega$
   for any $X\in T_x (\partial \Omega)$,
 where $\nabla \mathfrak{F}(x)=\sum_{j,l}\frac{\partial \mathfrak{F}}{\partial
  x_j}  g^{jl}\frac{\partial}{\partial x_l}$.
   In the local coordinates, the {\it Hessian} of $\mathfrak{F}(x)$ is
  \begin{eqnarray}  \nabla^{\Omega} \nabla  \mathfrak{F}=
  \left(\frac{\partial^2 \mathfrak{F} }{\partial x_i \partial x_j}
  -\frac{\partial \mathfrak{F}}{\partial x_k}
  \Gamma_{ij}^k\right),\end{eqnarray}
   and we have
   \begin{eqnarray}  \nabla^{\Omega} \nabla \mathfrak{F} (X, Y)=
   \langle \nabla_X^\Omega \nabla  \mathfrak{F}, Y\rangle, \quad \; X,Y\in T_x, \end{eqnarray}
   where $\Gamma_{ij}^k$
    is the {\it Christoffel symbols}.
   Therefore, $-\nabla^\Omega \nabla \mathfrak{F}$ has $n$
   eigenvalues at $x\in \partial \Omega$, one of which is $0$
   (corresponding to the eigenvector $\nabla \mathfrak{F}(x)$),
   and the others are the principal curvatures of $\partial \Omega$ at $x$.

\vskip 0.26 true cm

Let us consider the curvature of the boundary of a geodesic ball
$B_r(x_0)$ in ${\Bbb M}_k$. Since any two geodesic balls with the
same radius in ${\Bbb M}_k$ are isometric, their boundaries have the
same curvature. It is easy to check that the geodesic sphere of
radius $r$ with center at the origin has constant curvature $\tau_k
(r)$ (see [9, p.66]),
\begin{eqnarray}  \tau_k(r)=\left\{ \begin {array}{ll}  \sqrt{k}\cot \sqrt{k} \,r
  & \quad \mbox {if}\;\; k>0,\\
  \frac{1}{r}  & \quad \mbox{if}\;\;  k=0,\\
\sqrt{-k}\coth \sqrt{-k}\, r &  \quad \mbox{if}\;\;
 k<0.\end{array}\right.\end{eqnarray}

\vskip 0.28 true cm

Let $\Omega$ be a domain with $C^2$ boundary in either Euclidean
space ${\Bbb R}^n$, or the hyperboloid model ${\Bbb H}^n_\rho$, or
 the sphere ${\Bbb S}^n_\rho$. In the
    last case, $\Omega$ is required to lie in a hemisphere.
   Let $\Omega$ contain $x_0$, where $x_0$ is either the origin in Euclidean space ${\Bbb R}^n$, or
   the south pole $(0, \cdots, 0, -\rho)$ on the sphere ${\Bbb S}^n_\rho$,
   or the point
   $(0, \cdots, 0, \rho)$) on the hyperboloid model ${\Bbb
   H}^n_{\rho}$.
     We define the orthogonal projection $P_0$ from $\Omega$ to the Euclidean space
           $\{(x,0)\in {\Bbb R}^{n+1}\big|x\in {\Bbb R}^n\}$ by
    \begin{eqnarray} \;  P_0(y)= \left\{ \begin{array}{ll}
    (y_1, \cdots, y_{n}), \, \quad & \forall \; y=(y_1,\cdots, y_{n+1})\in \Omega\subset {\Bbb S}^n_{\rho},\\
    (y_1, \cdots, y_n), \, \quad & \forall \; y=(y_1, \cdots, y_n)\in \Omega \subset {\Bbb R}^n, \\
    (y_1,\cdots, y_n), \, \quad & \forall \; y=(y_1, \cdots, y_{n+1})\in \Omega \subset
    {\Bbb H}^n_\rho.\end{array} \right.
   \end{eqnarray}

  \noindent{Lemma 4.1.} \   {\it  Let  $\Omega$, $x_0$ and $P_0$ be as in the
  above description.
     Assume that $q\in \partial \Omega\cap \partial B_R
   (x_0)$, where $B_R(x_0)\subset \Omega$ is an open
   geodesic ball of ${\Bbb M}_k$ with geodesic radius $R>0$
   center at $x_0$.
   Denote by $\lambda_i(q)$ (respectively, ${\tilde
            \lambda}_i(P_0(q))$)
    the principal curvatures of $\partial \Omega$ at $q$ (respectively, $P_0(\partial \Omega)$
    at $P_0(q)$). Then
 \begin{eqnarray} \lambda_i(q) =\big({\tilde \lambda}_i (P_0(q))\big) h'_k(R),\end{eqnarray}
 where \begin {eqnarray} h'_k(r)= \left\{ \begin{array}{ll}
    \cos \sqrt{k}\, r, \quad & k>0,\\
    1,  \quad & k=0, \\
    \cosh \sqrt{-k} \,r, \quad & k<0.\end{array} \right.
   \end{eqnarray}}

\vskip 0.23 true cm

\noindent {\it Proof.}   \ \   It suffice to prove this lemma for
  spherical and hyperboloid model cases.

 \vskip 0.1 true cm

  (i) \ \  For spherical case,
  recall that ${\Bbb S}^n_\rho=\{y\in {\Bbb R}^{n+1}\big| \sqrt{y_1^2 +\cdots
 +y_{n+1}^2}=\rho\}$ of radius $\rho=1/\sqrt{k}$,
 centered at the origin in ${\Bbb R}^{n+1}$, with the induced
 Euclidean metric.
     Let $\{e_1, \cdots, e_{n-1}, \nu\}$ be a local orthonormal frame filed
     in a neighborhood of $q$
   such that $e_1, \cdots, e_{n-1}$ are the principal curvature vectors of
   $\partial \Omega$ and $\nu$ is the  exterior unit normal vector to
     the boundary $\partial \Omega$ of $\Omega$.
             Since $\langle \nu, e_j\rangle=0$, we get that
   $$0\equiv e_i \langle \nu, e_j\rangle = \langle \nabla_{e_i}
     \nu, e_j\rangle +\langle \nu, \nabla_{e_i} e_j\rangle$$
 for all $i,j=1, \cdots, n-1,\,$
 i.e.,
  \begin {eqnarray} I\!I(e_i, e_j):= \left\langle \nabla_{e_i} \nu, e_j\right\rangle
 =  - \left\langle \nu, \nabla_{e_i} e_j\right\rangle,\end{eqnarray}
 where $I\!I$ is the second fundamental form of $\partial
\Omega$, and the inner product $\langle \cdot, \cdot\rangle$ is
taken in the induced Euclidean metric. Similarly, we have
 \begin{eqnarray} {I\!I}({\tilde e}_i, {\tilde e}_j) = \left\langle \nabla_{{\tilde e}_i} {\tilde \nu},
 {\tilde e}_j\right\rangle
 =  - \left\langle {\tilde \nu}, \nabla_{{\tilde e}_i} {\tilde e}_j\right\rangle,\end{eqnarray}
  where  $\{{\tilde e}_1, \cdots, {\tilde e}_{n-1}, {\tilde \nu}\}$ is a local  orthonormal frame filed
     in a neighborhood of $P_0 (q)$ in Euclidean space $\{(x,0)\in {\Bbb R}^{n+1} \big|x\in {\Bbb R}^n\}$
      such that ${\tilde e}_1, \cdots, {\tilde e}_{n-1}$ are the principal curvature vectors of
   $P_0(\partial \Omega)$ and ${\tilde \nu}$ is the exterior unit normal vector to
     the boundary $P_0(\partial \Omega)$ of $P_0(\Omega)$.
 For any $y\in \Omega \subset {\Bbb S}^n_\rho$, it is obvious (see, for example, [8, p.62]) that
  $$y=\big((\frac{1}{\sqrt{k}}\sin \sqrt{k} r) \,\theta,
  \frac{1}{\sqrt{k}}\cos \sqrt{k}\, r),$$
  and hence  $$P_0 (y) = (\frac{1}{\sqrt{k}}\sin
 \sqrt{k} \, r)\,\theta,$$
  where $\theta\in {\Bbb S}^{n-1}$ and $r$ is the geodesic distance from the south pole $x_0$ to $y$.
   By our assumption,
it follows that $\nu(q)=\big((\cos \sqrt{k}\, R)\theta,
-\sin\sqrt{k}\, R\big)$.
 Thus, in the Euclidean space ${\Bbb R}^{n+1}$ we
   have  that $$e_i(q)={\tilde e}_i (P_0 (q)) \quad \; \mbox{for all } \, i=1, \cdots, n-1,$$  and
  $$\langle \nu(p), {\tilde \nu} (P_0 (q))\rangle =\langle \nu(q), (\theta, 0)\rangle =\cos \sqrt{k} \,
  R.$$  From this and (4.8)--(4.9), we get the corresponding part of (4.6) for $k>0$.

\vskip 0.24 true cm

   (ii) \ \   Recall that  $${\Bbb H}^n_\rho =\{y\in {\Bbb
R}^{n+1}\big|\langle y,
 y \rangle=-\rho,\;  y_{n+1}>0\}$$
  with the Riemannian metric induced from the Lorentzian metric
    $$\langle y,y\rangle = -y_{n+1}^2 +y_1^2+\cdots y_n^2,$$
 where $\rho=1/\sqrt{-k}$.
          Let $\{e_1, \cdots, e_{n-1}, \nu\}$ be a local orthonormal frame filed
     in a neighborhood of $q$ such that $e_1, \cdots, e_{n-1}$ are the principal curvature vectors of
   $\partial \Omega$
  and $\nu$ is the  exterior unit normal vector to
     the boundary $\partial \Omega$ of $\Omega$.
        Since $\langle \nu, e_j\rangle=0$, we get that
    \begin {eqnarray*} I\!I(e_i, e_j) = \left\langle \nabla_{e_i} \nu, e_j\right\rangle
 =  - \left\langle \nu, \nabla_{e_i} e_j\right\rangle,\end{eqnarray*}
 where $\langle \cdot, \cdot\rangle$ is taken in the
   Lorentzian metric.   Similarly, we have
 \begin{eqnarray*} {I\!I}({\tilde e}_i, {\tilde e}_j) =
 \left\langle \nabla_{{\tilde e}_i} {\tilde \nu}, {\tilde e}_j\right\rangle
 =  - \left\langle {\tilde \nu}, \nabla_{{\tilde e}_i} {\tilde e}_j\right\rangle,\end{eqnarray*}
  where  $\{{\tilde e}_1, \cdots, {\tilde e}_{n-1}, {\tilde \nu}\}$ is a local  orthonormal frame filed
     in a neighborhood of $P_0 (q)$ in Euclidean space $\{(x,0)\in {\Bbb R}^{n+1}\big|x\in {\Bbb
     R}^n\}$
      such that ${\tilde e}_1, \cdots, {\tilde e}_{n-1}$ are the principal curvature vectors of
   $P_0(\partial \Omega)$ and ${\tilde \nu}$ is the exterior unit normal vector to
     the boundary $P_0(\partial \Omega)$ of $P_0(\Omega)$.
  Note that  for any $y\in \Omega \subset {\Bbb H}^n_\rho$, one has (see, for example, [10,
 p.22]) that
 \begin{eqnarray*} y=\left(\big(\frac{1}{\sqrt{-k}}\sinh \sqrt{-k}\,r\big)\theta,
  \frac{1}{\sqrt{-k}}\cosh \sqrt{-k}\,r\right),\end{eqnarray*}
and hence
  $$P_0 (y) =\big(\frac{1}{\sqrt{-k}}\sinh \sqrt{-k}\,r\big) \theta \quad \;\; \mbox{for any} \; \;  y\in
  \Omega,$$
   where $\theta\in {\Bbb S}^{n-1}$ and $r$ is the geodesic distance from the point $x_0=(0,\cdots, 0, \rho)$ to
   $y$.
   By the assumption, in the Euclidean space ${\Bbb R}^{n+1}$ with Lorentzian metric, we
   then have that
   $$e_i(q)={\tilde e}_i (P_0 (q))\quad \; \mbox{for all } \; i=1, \cdots, n-1,$$ and
     $$\nu(q)=\big((\cosh \sqrt{-k}\,R)\theta, \sinh
     \sqrt{-k}\,R\big),$$ which implies
       $$\langle \nu(q), {\tilde \nu} (P_0 (q))\rangle =
\langle \big((\cosh \sqrt{-k}\,R)\theta, \sinh
     \sqrt{-k}\,R\big), (\theta, 0)\rangle= \cosh
  \sqrt{-k} \, R.$$
   Therefore, we obtain the corresponding part of (4.6).

\vskip 0.26 true cm

  \noindent{Theorem 4.2.} \ \  {\it Let $\Omega$
  be a domain with $C^2$ boundary in the $n$-dimensional space
  ${\Bbb M}_k$ of constant curvature $k$, $n\ge 2$, and let $\lambda_1, \cdots, \lambda_{n-1}$
  denote the principal curvatures of $\partial \Omega$.
   Assume that $B_R(x_0)\subset \Omega$ is an open geodesic ball with radius $R>0$
   center at $x_0$ and suppose that the set
   $\partial \Omega\cap \partial B_R (x_0)$
  is made of a finite number  of points $p_1, \cdots, p_l$
  such that $\lambda_j(p_m)<\tau_k(R)$ for every $j=1, \cdots, n-1$ and
  every  $m=1, \cdots, l$,  where $\tau_k(R)$ is as in (4.4).
 Let $W=W(s, x)$ be the solution to problem
 \begin{eqnarray}  &&\sum_{i=1}^n \left(\frac{1+k|x|^2}{2}\right)^2 \, \frac{\partial^2 W}
   {\partial x_i^2} -s W=0 \quad \; \mbox{in}\;\, \Omega, \\
   && \; W=1 \quad \;\; \quad \quad \;\quad \mbox{on}\;\; \partial
 \Omega.\end{eqnarray}
   Then, the following formula holds for every function $\phi$
   continuous on ${\Bbb M}_k$:
   \begin{eqnarray} && \lim_{s\to +\infty} s^{\frac{n-1}{4}} \int_{\partial B_R(x_0)}
   \phi(x)W(s, x) \, dA_x\\
   && \quad \;\; =(2\pi)^{\frac{n-1}{2}} \sum_{m=1}^l
   \phi(p_m) \left[\frac{1}{(h'_k(R))^{n-1}}\prod_{j=1}^{n-1} \bigg(\tau_k(R) -\lambda_j
   (p_m)\bigg)\right]^{-1/2},\nonumber \end{eqnarray}
 where $h'_k(r)$ is as in (4.7)}.

\vskip 0.23 true cm

 \noindent {\it Proof.} \ \    Let $p_m\in \{p_1, \cdots, p_l\}$; by applying a
 partition of unity, we can suppose that supp$\,\phi$ does not contain any
 $p_i$ different from $p_m$.

 Since there exists an isometry $\Phi$ that maps ${\Bbb M}_k$
 onto itself such that $\Phi x_0=0$
 and the equation (4.10) is invariant under
  the isometry map $\Phi$,
   we may assume that $x_0=0$
  and use the spherical coordinates about the point $x_0=0$.
   As in (3.9),  we have
\begin{eqnarray} && \int_{\partial B_R (0)} \phi(x) e^{-\sqrt{s}\, \mathfrak{F}(x)}
   dA_x  \\
   && \;\;\qquad =
\int_{{\Bbb S}^{n-1}} \big(h_k(R)\big)^{n-1} \phi(\exp_0 (R\theta))
e^{-\sqrt{s}\,
  \mathfrak{F}(\exp_0(R\theta))}
   d\Theta(\theta)\nonumber \\
     && \;\;\qquad =\int_{{\Bbb S}^{n-1}_{h_k(R)}} \phi\big(P_0^{-1} ({\tilde x})\big) e^{-\sqrt{s}\,
 \mathfrak{F}(P_0^{-1}({\tilde x}))}
   d\Theta({\tilde x})\nonumber,\end{eqnarray}
 where ${\Bbb S}_{h_k(R)}^{n-1}$ is the sphere of radius $h_k(R)$
 with center at the origin in Euclidean space ${\Bbb R}^n$,
   $P_0^{-1}$ is the inverse of $P_0$. Here $P_0$ is the
   orthogonal projection from $\Omega$ to the
    Euclidean space ${\Bbb R}^n$ as before
     (Note that in order to say $P_0$, we must regard ${\Bbb M}_k$ as either
    ${\Bbb R}^n$, or the sphere ${\Bbb S}^n_\rho$, or the hyperboloid model ${\Bbb
    H}^n_\rho$).
   For convenience, we denote by ${\tilde x}$ the point $P_0(x)$ for any $x\in
   \Omega$.
       Also, we can suppose that
$P_0^{-1}(supp \;\phi)$ does not contain the point $\,-P^{-1}_0
(p_m)$.
  As in [23, p.938], let ${\Bbb
R}^{n-1} \ni \eta =(\eta_1, \cdots, \eta_{n-1})\mapsto {\tilde
x}(\eta)\in {\Bbb S}^{n-1}_{h_k(R)}$ be the stereographic projection
from the point
  $-{\tilde p}_m$ onto the tangent space to
  ${\Bbb S}^{n-1}_{h_k(R)}$ at ${\tilde p}_m$.
  More precisely, take an orthogonal basis $\xi^1, \cdots, \xi^n$ of
  ${\Bbb R}^n$ with $\xi^n= -\frac{{\tilde p}_m}{h_k(R)}$, and put
$$ {\tilde x}(\eta)=\frac{2h_k(R) |\eta|^2}{(2h_k(R))^2 +|\eta|^2} \xi^n+
\frac{(2h_k(R))^2}{(2h_k(R))^2+|\eta|^2} \sum_{j=1}^{n-1} \eta_j
\xi^j +{\tilde p}_m.$$  Thus, we have
\begin{eqnarray} \;\; \qquad  \int_{\partial B_R (0)} \phi(x) e^{-\sqrt{s}\, \mathfrak{F}(x)}
   dA_x  =
 \int_{{\Bbb R}^{n-1}}  \phi(P_0^{-1}({\tilde x}(\eta))) e^{-\sqrt{s}
\,
    \mathfrak{F} (P_0^{-1}({\tilde x}(\eta)))}  J(\eta)d\eta,\end{eqnarray}
   where $$ J(\eta) := \sqrt{\mbox{det}\, \left(\frac{\partial
   {\tilde x}(\eta)}{\partial \eta_i}\cdot \frac{\partial
   {\tilde x}(\eta)}{\partial \eta_j}\right)}
   =\left(\frac{(2h_k(R))^2}{(2h_k (R))^2+|\eta|^2}\right)^{n-1}.$$
 Set $\mathfrak{F}^*(\eta)= (\mathfrak{F}\circ P_0^{-1})({\tilde x}(\eta))$.
 Then $\mathfrak{F}^*(0)=0$, and
 ${\tilde \nabla} \mathfrak{F}^*(0)=0$ and ${\tilde \nabla}^2 {\mathfrak{F}}^* (0)$ is
 positive definite,
where ${\tilde \nabla}$ and ${\tilde \nabla}^2$ is in the sense of
Euclidean metric. In fact,
  differentiating
${\mathfrak{F}}^*(\eta)$ twice yields (cf. [23, p.938]):
  \begin{eqnarray}  \quad \frac{\partial^2 {\mathfrak{F}}^*}{\partial \eta_i\partial
  \eta_j}(\eta)&=&\frac{\partial {\tilde x}}{\partial
  \eta_i}(\eta)\cdot \left(\big({\tilde \nabla}^2 (\mathfrak{F}\circ P_0^{-1}) ({\tilde
  x}(\eta))\big)
  \frac{\partial  {\tilde x}}{\partial \eta_j}(\eta)\right)\\
  & &+ \big({\tilde \nabla} (\mathfrak{F}\circ P_0^{-1})({\tilde x}(\eta))\big)
   \cdot \frac{\partial^2 {\tilde x}}{\partial
  \eta_i\partial \eta_j}(\eta), \quad \, i, j=1, \cdots,
  n-1,\nonumber\end{eqnarray}
  for every $\eta\in {\Bbb R}^{n-1}$, where the dot denotes scaler
  product of vectors in ${\Bbb R}^n$.
        It follows from ${\tilde x}(\eta)\in {\Bbb S}^{n-1}_{h_k(R)}$ for every
    $\eta\in {\Bbb R}^{n-1}$ that
    \begin{eqnarray*} &&\frac{\partial {\tilde x}}{\partial \eta_i}
 (\eta)\cdot \big(({\tilde x} (\eta))-0\big)=0,\quad \; i=1,\cdots, n-1,\\
   &&\frac{\partial^2 {\tilde x}}{\partial \eta_i \partial
   \eta_j}(\eta)\cdot\big(({\tilde x}(\eta))-0\big)+\frac{\partial {\tilde x}}{\partial
   \eta_i}(\eta)\cdot \frac{\partial {\tilde x}}{\partial
   \eta_j}(\eta)=0, \quad \, i,j=1,\cdots, n-1,\end{eqnarray*}
    for all $\eta\in {\Bbb R}^{n-1}$.
Clearly, as in [23, p.939] we have
   \begin{eqnarray*}
   -\big({\tilde \nabla} (\mathfrak{F}\circ P_0^{-1})\big) ({\tilde p}_m)
  =\frac{({\tilde x}(0)-0)}{h_k(R)}.\end{eqnarray*}
   Thus  \begin{eqnarray} &&\big(({\tilde \nabla} (\mathfrak{F}\circ P_0^{-1}))({\tilde
  p}_m)\big)
  \cdot \frac{\partial {\tilde
  x}}{\partial \eta_i}(0)= 0, \quad i=1, \cdots, n-1,\\
   &&\big(({\tilde \nabla} (\mathfrak{F}\circ P_0^{-1}))({\tilde p}_m)\big)
   \cdot \frac{\partial^2 {\tilde x}}{\partial
   \eta_i\partial \eta_j}(0)=
   \left(\frac{1}{h_k(R)}\right) \frac{\partial {\tilde
   x}}{\partial \eta_i}(0) \cdot \frac{\partial {\tilde x}}{\partial
   \eta_j}(0),\\
   && \quad  \; \; i, j=1, \cdots, n-1. \nonumber \end{eqnarray}
    We find from this and (4.15) that
    \begin{eqnarray*}  && {\tilde \nabla} {\mathfrak{F}}^*(0)=0,\\
   &&\frac{\partial^2 {\mathfrak{F}}^*}{\partial \eta_i\partial \eta_j} (0)=
    \frac{1}{h'_k(R)}\left\{\frac{\partial {\tilde x}}{\partial \eta_i}(0) \cdot \left[
  \left((h'_k(R)) \big(({\tilde \nabla}^2 (\mathfrak{F}\circ
   P_0^{-1})) ({\tilde p_m})\big)+\frac{h'_k(R)}{h_k(R)}\, I \right)
 \frac{\partial {\tilde x}}{\partial \eta_j}(0)\right]\right\}, \\ &&
 \quad \;  \quad \;  i,j=1,\cdots, n-1,\end{eqnarray*}
  where $I$ is the $(n-1)\times (n-1)$ identity matrix.
    Since \begin{eqnarray*} \frac{\partial {\tilde x}(\eta)}{\partial \eta_i}
\; \frac{\partial {\tilde x}(\eta)}{\partial \eta_j}=
\left(\frac{(2h_k(R))^2}{(2h_k (R))^2+|\eta|^2}\right)\delta_{ij},
\quad \; i,j= 1,\cdots, n-1,\end{eqnarray*}
 we see that the vectors
$\frac{\partial {\tilde x}}{\partial
 \eta_i} (0)$,
 $i=1, \cdots, n-1$, make an orthogonal basis of the tangent space
 $T_{{\tilde p}_m} (P_0(\partial \Omega))
 =T_{{\tilde p}_m}\big(P_0(\partial B_R(0))\big)$.
 Thus
\begin{eqnarray*} && {\tilde \nabla} {\mathfrak{F}}^*(0)=
   0,\\
  &&  det\, {\tilde \nabla}^2 {\mathfrak{F}}^*(0)=
   \frac{1}{(h'_k (R))^{n-1}} \, det \left(h'_k(R)\big(({{\tilde \nabla}}^2 (\mathfrak{F}\circ
   P_0^{-1}))
  ({\tilde p}_m)\big)+ \tau_k(R) I\right).\end{eqnarray*}
   Let $0$ and $\{{\tilde \lambda}_j ({\tilde p}_m)\}_{j=1}^{n-1}$ be the
   eigenvalues of matrix $-\big({\tilde \nabla}^2 (\mathfrak{F}\circ
   P_0^{-1})\big)({\tilde p}_m)$.  Clearly, $\{{\tilde
   \lambda}_j({\tilde p}_m)\}_{j=1}^{n-1}$ are the principal curvatures of the
   boundary $P_0(\partial \Omega)$ of $P_0(\Omega)$ at ${\tilde
   p}_m$.
     It follows from Lemma 4.1 that, under map $P_0$,
     $\{{\tilde \lambda}_j ({\tilde p}_m)\}_{j=1}^{n-1}$ and the principal curvature
     $\{\lambda_j(p_m)\}_{j=1}^{n-1}$
  of the boundary $\partial \Omega$ of $\Omega$ at $p_m$
  have the following formula:
 \begin{eqnarray*} \lambda_j(p_m)=({\tilde \lambda}_j({\tilde p_m}))\,
 h'_k(R), \quad \, j=1,\cdots, n-1, \end{eqnarray*}
where $R=dist(0, p_m)$, and $h'_k(R)$ is as in (4.7).
    Note that  $(\mathfrak{F}\circ P_0^{-1}) {\tilde x}=\mathfrak{F}(x)$
    for every $x\in \bar\Omega$.
Since the eigenvalues of matrix $-(\nabla^2 \mathfrak{F}(p_m)
 +\tau_k(R)I)$ are $0$ and
  $(\tau_k(R)-\lambda_j(p_m))$,  $j=1, \cdots, n-1$, where
  $\tau_k(R)$ is the constant curvature of the geodesic sphere
  $\partial B_R(0)$ in ${\Bbb M}_k$,
   it follows that
   \begin{eqnarray*}&&det \, {\tilde \nabla}^2
   {\mathfrak{F}}^*(0)=
 \frac{1}{(h'_k(R))^{n-1}}\nonumber \\
     &&\times \left|\begin{array}{cccc}
  \tau_k(R)- \big({\tilde \lambda}_1({\tilde p}_m)\big) h'_k(R) & 0 & \cdots & 0\\
      0 &  \tau_k(R)-\big({\tilde \lambda}_2({\tilde p}_m)\big) h'_k(R)  &  \cdots &  0\\
      \vdots & \vdots & \ddots &  \vdots \\
      0 & 0 &\cdots &  \tau_k(R)-\big({\tilde \lambda}_{n-1}({\tilde p}_m)\big)
      h'_k(R)
  \end{array} \right|\nonumber \\
  &&{}\\
 &&=\frac{1}{(h'_k(R))^{n-1}}\left|\begin{array}{cccc}
   \tau_k(R)-\lambda_1(p_m)  & 0 & \cdots & 0\\
      0 &  \tau_k(R)-\lambda_2 (p_m)   &  \cdots &  0\\
      \vdots & \vdots & \ddots &  \vdots \\
      0 & 0 &\cdots &  \tau_k(R)-\lambda_{n-1}(p_m)
  \end{array} \right| \nonumber\\
   &&= \frac{1}{(h'_k(R))^{n-1}}
   \prod_{j=1}^{n-1}
  \left[\tau_k(R)-\lambda_j(p_m)\right], \end{eqnarray*}
  i.e.,
   \begin{eqnarray} \mbox{det}\, {\tilde \nabla}^2 {\mathfrak{F}}^*(0)=\frac{1}{(h'_k(R))^{n-1}}
           \, det \big(\nabla^2 \mathfrak{F}(p_m)+ \tau_k(R)
   I\big). \end{eqnarray}
  (Note that $0$ and the principal curvatures
$\{\lambda_j(p_m)\}_{j=1}^{n-1}$ of $\partial \Omega$ at $p_m$
 are the eigenvalues of the Hessian matrix [16,
  p.139]
   $$ -\nabla^2 \mathfrak{F}(p_m):=-
 \left(\frac{\partial^2 \mathfrak{F}}{\partial x_i \partial x_j}
 -\frac{\partial \mathfrak{F}}{\partial x_q} \Gamma_{ij}^q\right)\bigg|_{p_m}\,).$$

Since supp $\phi$ does not contain any $p_i$ different from $p_m$,
 we may assume that ${\mathfrak{F}}^* (\eta)>0$ if $\eta\ne 0$.
 Hence by  Laplace's method (see [6, p.71]), or by the stationary phase method
 (see [11, p. 208-217] or [23] for example),
\begin{eqnarray}  &&  \lim_{s\to +\infty} s^{\frac{n-1}{4}}
\int_{{\Bbb R}^{n-1}} {\phi}(P_0^{-1} ({\tilde x}(\eta)))
e^{-\sqrt{s}\,
 \big((\mathfrak{F}\circ P_0^{-1}) ({\tilde x}(\eta))\big)}
   J(\eta) d\eta \\
   && \quad \;\;  =  (2\pi)^{\frac{n-1}{2}} {\phi}(p_m)
   J(0) \big(det{\tilde \nabla}^2 {\mathfrak{F}}^* (0)
   \big)^{-\frac{1}{2}}.\nonumber\end{eqnarray}
  From $ J(0)=1$, (4.14), (4.18) and (4.19), we get
\begin{eqnarray}  &&  \lim_{s\to +\infty} s^{\frac{n-1}{4}}
\int_{\partial B_R (0)} \phi(x) e^{-\sqrt{s}\, \mathfrak{F}(x)}
   dA_x \\
   && \quad \;\;  =  (2\pi)^{\frac{n-1}{2}} \sum_{m=1}^l \phi(p_m)
  \bigg(\frac{1}{(h'_k(R))^{n-1}}\prod_{j=1}^{n-1} \big[\tau_k (R) -\lambda_j
  (p_m)\big]\bigg)^{-1/2}.\nonumber\end{eqnarray}

Finally, we prove formula (4.11).  It is sufficient to prove it
 for any nonnegative function  $\phi$ (see [23, p.940]).
 From Lemma 3.3, one has that for all $s\ge s_\epsilon$ and any
nonnegative $\phi$,
$$ \int_{\partial B_R(x_0)} \phi(x) W_\epsilon^- (s, x) dA_x \le
 \int_{\partial B_R(x_0)} \phi(x) W(s, x) dA_x \le
\int_{\partial B_R(x_0)} \phi(x) W_\epsilon^+ (s, x) dA_x.$$
  Therefore, (4.19) and the definition (3.20) implies that
\begin{eqnarray*}  &&
\bigg(\frac{2\pi}{\sqrt{1+\epsilon}}\bigg)^{\frac{n-1}{2}}
   \sum_{m=1}^l \phi(p_m)
  \bigg(\frac{1}{(h'_k(R))^{n-1}}
  \prod_{j=1}^{n-1} \big[\tau_k (R) -\lambda_j
  (p_m)\big]\bigg)^{-1/2}\\
  &&  \qquad \;\; \le
 {\lim \inf}_{s\to +\infty} s^{\frac{n-1}{4}}
\int_{\partial B_R (x_0)} \phi(x) W(s,x)dA_x \\
 && \qquad \; \; \le {\lim \sup}_{s\to
+\infty} s^{\frac{n-1}{4}}
 \int_{\partial B_R (x_0)} \phi(x) W(s,x)dA_x \\
   && \qquad \;\;  \le  \bigg(\frac{2\pi}{\sqrt{1-\epsilon}}\bigg)^{\frac{n-1}{2}}
   \sum_{m=1}^l \phi(p_m)
  \bigg(  \frac{1}{(h'_k(R))^{n-1}}
  \prod_{j=1}^{n-1} \big[\tau_k (R) -\lambda_j
  (p_m)\big]\bigg)^{-1/2}.\end{eqnarray*}
 for every $\epsilon>0$. By letting $\epsilon$ tend to $0$, we
 get (4.11) and the proof is completed.  $\;\, \square$

 \vskip 1.32 true cm

 \vskip 1.32 true cm

\section{Proof of main results}

\vskip 0.43 true cm

In this section, we shall prove the analyticity of the boundary
$\partial \Omega$ and the main theorem.  A domain $\Omega$ is said
to satisfy the exterior geodesic sphere condition if for every $y\in
\partial \Omega$ there exists a geodesic ball $B_r(z)$ such that
$\overline{B_r(z)}\cap \bar \Omega=y$. A domain $D$ satisfies the
interior geodesic cone condition if for every $x\in
\partial D$ there exists a finite geodesic spherical cone $K_x$
with vertex $x$ such that $K_x \subset \bar D$ and $\bar K_x \cap
\partial D=\{x\}$.

\vskip 0.23 true cm

  \noindent{Lemma 5.1.} \ \  {\it Let $\Omega$ be a
  bounded domain in $n$-dimensional space ${\Bbb
  M}_k$ of constant curvature $k$. Let $\Omega$ satisfy the
  exterior geodesic sphere condition and suppose that
  $D$ is a domain satisfying the interior geodesic cone
  condition and such that $\bar D
  \subset \Omega$. Assume that the solution $u=u(t, x)$ of problem
  (1.1)--(1.3) satisfies condition (1.4). Let $R$ be the positive
  constant given by
  \begin{eqnarray}  R=\lim_{s\to +\infty} \left[-\frac{1}{\sqrt{s}} \log
  A(s)\right],\end{eqnarray}
  where \begin{eqnarray}   W(s, x)=s\int_0^{+\infty} a(t) e^{-s\,t} dt:=
  A(s), \quad \, x\in \partial D.\end{eqnarray}
 Then the following assertions hold:

 (i)  \ \  for every $x\in \partial D$,  $\mathfrak{F}(x)=R$, where $\mathfrak{F}$ is
 defined by (2.10);

 (ii)  \ \  $\partial D$ is real analytic;

 (iii) \ \  $\partial \Omega$  is real analytic and $\partial \Omega= \{x\in {\Bbb
 M}_k \,\big| \, \mbox{dist}\, (x,\partial D)=R\}$;}

(iv) \ \   Let $\lambda_j(y)$, $j=1, \cdots, n-1$ denote the
$j^{th}$ principal curvature at $y\in \partial \Omega$  of the
 real analytic surface $\partial \Omega$; then $\lambda_j(y)<\tau_k(R)$,
$j=1, \cdots, n-1$, for every $y\in
\partial \Omega$, where $\tau_k(R)$ is given by (4.4).

\vskip 0.24 true cm

\noindent{Proof.} \ \  (i) \ \  Let
  \begin{eqnarray} W(s,x)=s \int_0^{+\infty} u(t, x)e^{-s\,t}dt,
 \quad \, s>0.\end{eqnarray}
  Then $W(s,x)$ satisfies elliptic boundary value problem
  (3.18)--(3.19). Applying Lemma 2.1 (Varadhan's theorem), we have
  \begin{eqnarray} \lim_{s\to +\infty} \left(-\frac{1}{\sqrt{s}}
  \log W(s, x)\right)= d(x, \partial \Omega)=\mathfrak{F}(x).\end{eqnarray}
 Since $u$ satisfies (1.4), it follows that for fixed $s>0$,
  $A(s)$ is constant on $\partial D$. Therefore,  $\mathfrak{F}(x)=R$ for every
   $x\in \partial D$.

  (ii) \ \ It follows from the interior regularity of parabolic equations (see \cite{Mo2})
  that $u(t, x)$ is real analytic on any compact subdomain in $\Omega$.
   By the implicit function theorem for real analytic function
  (see, for example, [15, p.69]), it suffices to prove that, for every point $x\in
  \partial D$, there exists a time $t^*>0$ such that $\nabla u(t^*,
  x)\ne 0$.  Suppose by contradiction that there exists a  point $x_0\in
     \partial D$ such that $\nabla u(t,x_0)=0$ for every $t>0$.
     It follows from Lemma 3.2 (ii) that
     \begin{eqnarray} \int_{\partial B_R(x_0)} \exp_{x_0}^{-1}
     (x-x_0) u(t, x) dA_x=0 \quad \; \mbox{for every } t>0.
     \end{eqnarray}
 We may put $x_0=0$ by an isometry of ${\Bbb M}_k$. Thus, (5.5) is  equivalent to
  \begin{eqnarray*} \int_{{\Bbb S}^{n-1}_{h_k(R)}} \xi  u(t, \xi)
  d\Theta(\xi)=0 \quad \; \mbox{for every } t>0,
     \end{eqnarray*}
and hence
 \begin{eqnarray} \int_{{\Bbb S}^{n-1}_{h_k(R)}} \xi  W(s, \xi)
  d\Theta(\xi)=0 \quad \; \mbox{for every } s>0.
     \end{eqnarray}
 On the other hand, as in [23, p.941-942] one can show that
 \begin{eqnarray*} \int_{{\Bbb S}^{n-1}_{h_k(R)}} \xi  W(s, \xi)
  d\Theta(\xi)>0
     \end{eqnarray*}
for $s>0$ sufficiently large. This is a contradiction.

 (iii)  \ \   Set $\Gamma=\{y\in {\Bbb M}_k\big| d(y, D)=R\}$.
  For each $y\in \Gamma$ there exists a point $x\in \partial D$ such
  that $d(y, D)=d(y, x)$.
   Let $\gamma(r)$ be a geodesic starting from $x$ and ending at $y$.
  We claim that $\dot{\gamma}(0):=\frac{d\gamma(r)}{dr}(0)$ is orthogonal to
  the tangent space of $\partial D$.
     In fact,  let $\zeta (\mu)$ be a smooth curve
   in $\partial D$ with $\zeta (0)=x$. For each $\zeta (\mu)$, let $\gamma_\mu (r)$
    be the geodesic starting from $\zeta (\mu)$ and ending at $y$, and
    let $L(\mu)$ be the length of the geodesic $\gamma_\mu (r)$ between $\zeta
(\mu)$ and $y$. Then $$L(\mu)=\int_0^{L(\mu)} \langle
  {\dot{\gamma}}_\mu(r), {\dot{\gamma}}_\mu(r)\rangle^{1/2} dr.$$
     It is easy to check
    that $L(\mu)$ has the following variational formula (cf. [8, p.67]):
 \begin{eqnarray*} \frac{dL(\mu)}{d\mu}\big|_{\mu=0} && =
  \left[\langle
  {\dot{\gamma}}_\mu(L(\mu)),
  {\dot{\gamma}}_\mu(L(\mu))\rangle^{1/2}\right]\left(\frac{dL}{d\mu}(0)\right)\\
  && \;\;  +\left[\langle {\dot{\gamma}}_0(r),
     \frac{\partial\gamma_\mu (r)}{\partial \mu}\big|_{\mu=0}\rangle
      \bigg|_{0}^{L(0)} -\int_0^{L(0)} \langle {\ddot{\gamma}}_0 (r),
     \frac{\partial \gamma_\mu (r)}{\partial \mu}\big|_{\mu=0}\rangle dr\right]\\
    &&= \langle{\dot{\gamma}}_0(L(0)), \frac{\partial \gamma_\mu (L(0))}
    {\partial \mu}\big|_{\mu=0}\rangle -
  \langle {\dot{\gamma}}_0 (0),
  \frac{\partial \gamma_\mu(0)}{\partial \mu}\big|_{\mu=0}\rangle.
  \end{eqnarray*}
  Here we have used the fact that
  $\frac{dL(\mu)}{d\mu}\big|_{\mu=0}=0$
  and $\ddot{\gamma}(r)=0\,$.
  From $\gamma_\mu (L(\mu))=y$ for all $\mu\ge 0$,
  we get $ \frac{\partial \gamma_\mu (L(0))}{\partial
  \mu}\big|_{\mu=0} = \frac{\partial \gamma_\mu (L(\mu))}{\partial
  \mu}\big|_{\mu=0}=0$,
  and hence
 $\langle {\dot{\gamma}}_0 (0), \frac{\partial
  \gamma_\mu (0)}{\partial \mu}\big|_{\mu=0}\rangle=
 \langle{\dot{\gamma}}_0 (0), \frac{d\zeta(\mu)}{d\mu}\big|_{\mu=0}\rangle=
0$. The claim is proved.

For the above $x\in \partial D$ (i.e., $d(x, y)=d(y, D)$), there
exists a unique
   point $y'\in \partial \Omega$ such that
   $\overline{B_R(x)}\cap \partial \Omega =\{y'\}$ (Indeed, if
   $y''\in \overline{B_R(x)}\cap \partial \Omega$ and $y'\ne
    y''$, then the geodesic $\beta_1(r)$ (connecting $x$ and $y'$)
    and the geodesic $\beta_2(r)$ (connecting $x$ and $y''$) have the
   same initial point $x$ and same direction at $x$. This is a
   contradiction). Since $\dot{\gamma} (0)$ is orthogonal
   to the tangent space of $\partial D$, it follows that $y=y'\in\partial \Omega$,
   and hence $\Gamma\subset \partial \Omega$.
  By the definition of $\Gamma$, we immediately derive that $\Gamma$
  is an analytic hypersurface diffeomorphic to $\partial D$.
  Therefore $\Gamma=\partial \Omega$, otherwise, $\partial \Omega$
  can't satisfy the exterior geodesic sphere condition.

    (iv) \ \   For any point $y\in \partial \Omega$, there exists a
    unique $x\in \partial D$ such that  $\overline{B_R(y)}\cap \bar
    D=\{x\}$.  Since $\partial D$ is real analytic, there exists
    a geodesic ball $B_r(z)\subset D$  such that
 $\overline{B_r(z)}\cap \partial D=\{x\}$.
 Thus,  $$\mathfrak{F}(z)=r+R \quad \mbox{and}\;\;
    \overline{B_{r+R}(z)}\cap \partial \Omega=\{y\},$$
    which implies
    $$\lambda_j (y)\le  \tau_k(R+r),\quad \; j=1, \cdots, n-1.$$
It is obvious by (4.4) that $\tau_k (R+r)< \tau_k (R)$. This
completes the proof of (iv).  $\;\; \square$

\vskip 0.46 true cm

 \noindent {\it Proof of Theorem 1.1.} \ \ From
Lemma 5.1, we
 see that $\partial \Omega$ and $\partial D$ are analytic.
  Let $p_1$  and $p_2$ be two distinct points in $\partial \Omega$.
  Then $\nabla \mathfrak{F}(p_i)$, is the unit interior normal vector
  of  $\partial \Omega$ at $p_i$, $i=1,2$.
    Let $\gamma_i (r)$  be the geodesic satisfying $\gamma_i(0)=p_i$ and
    ${\dot{\gamma}}_{i}(0)= \nabla \mathfrak{F}(p_i)$, $i=1,2$. It follows from Lemma 5.1
     that  $\gamma_i(R)\in \partial D$ and $\gamma_1 (R)\ne \gamma_2(R)$.
   Let us denote $\gamma_i(R)$ by $P_i$, and by $\Phi_i$ the isometric map of ${\Bbb M}_k$
   satisfying  $\Phi_i 0=P_i$, $i=1,2$.  Then for $x\in B_R(0)$,
   define the function $v(t, x)$ by
    \begin{eqnarray} v(t,x)= u(t, \Phi_1 x)- u(t, \Phi_2 x). \end{eqnarray}
  Lemma 3.1 implies  that $v(t, x)$ satisfies equation (1.1) in
  $(0,+\infty)\times B_R(0)$. By (1.4), we have
  $$v(t, 0)= u(t, P_1)- u(t,P_2)=0 \quad \, \mbox{for all} \; \,
  t>0.$$
 It follows from Lemma 3.2 (i) that
 $$ \int_{\partial B_R (0)} v (t, x) dA_x =0 \quad \;\; \mbox{for
 all } t>0,$$
 and hence
   $$ \int_{\partial B_R(P_1)} u(t, x)dA_x =  \int_{\partial B_R(P_2)} u(t,
   x)dA_x \quad \; \mbox{for all } t>0.$$
  Thus, by the definition of (5.3), we obtain
 \begin{eqnarray} \int_{\partial B_R(P_1)} W(s, x)dA_x =  \int_{\partial B_R(P_2)} W(s,
   x)dA_x \quad \; \mbox{for all } s>0.\end{eqnarray}
  Multiplying both sides of (5.8) by $s^{\frac{n-1}{4}}$ and letting $s\to
  +\infty$,
 by (4.12) of Theorem 4.2 (with $\phi\equiv 1)$ we get
   \begin{eqnarray*} \prod_{j=1}^{n-1} \left[ \tau_k(R)
   -\lambda_j(p_1)\right]=\prod_{j=1}^{n-1} \left[ \tau_k(R)
   -\lambda_j(p_2)\right], \end{eqnarray*}
which implies
\begin{eqnarray} \prod_{j=1}^{n-1} \left[ \tau_k (R)
   -\lambda_j(x)\right]=\mbox{constant},
   \quad \; \mbox{for every } x\in \partial \Omega. \end{eqnarray}
  Let us put
   \begin{eqnarray} F= F(\beta_1, \cdots, \beta_{n-1})=-\prod_{j=1}^{n-1} \left[ \tau_k(R)
   -\beta_j\right].\end{eqnarray}
 Clearly,  $F$ is of class $C^1$, and $F(\lambda_1, \cdots, \lambda_{n-1})=const$ on $\partial
 \Omega$.
  Since $\lambda_j(x)<\tau_k(R)$ on $\partial \Omega$, $\, j=1,\cdots, n-1$, we have
  $$const > \frac{\partial F(\beta_1, \cdots, \beta_{n-1})}{\partial \beta_j}
   >const>0 \quad \; (j=1, \cdots, n-1),$$
 at least on $\partial \Omega$, i.e., for $\beta_j=\lambda_j$ $\, (i=1,
 \cdots, n-1)$.
  It follows from Lemma 2.2 (Alexandrov's theorem) that
  $\partial \Omega$ must be a geodesic sphere in ${\Bbb M}_k$. $\;\; \square$

\vskip 0.25 true cm

 For the wave equations and the Schr\"{o}dinger
equations, we have the following:

\vskip 0.22 true cm

 \noindent  {\bf Theorem 5.2.} \ \  {\it  Let
$\Omega$ be a bounded domain in the $n$-dimensional space ${\Bbb
M}_k$ of constant curvature $k$ with the metric
$g_{ij}=\frac{4\delta_{ij}}{(1+k|x|^2)^2}\,$  (in case of $k>0$,
$\Omega$ is required to lie in a hemisphere), $\,n\ge 2$. Let
$\Omega$ satisfy the exterior geodesic sphere condition and assume
that $D$ is a domain, with boundary $\partial D$, satisfying the
interior geodesic cone condition, and such that $\bar D\subset
\Omega$.

   Suppose $v$ satisfies the following wave equation (5.11) or
   Schr\"{o}dinger's equation (5.12):
    \begin{eqnarray}  \left\{ \begin {array}{ll}  \frac{\partial^2 v}{\partial
  t^2}= \sum_{i=1}^n \frac{(1+k|x|^2)^2}{4} \,\,\frac{\partial^2 v}{\partial x_i^2}
  \quad \;\; &\mbox{in $\,(0,+\infty)\times \Omega$}\\
   v=1 \quad \,  &\mbox{on $\,(0,
  +\infty)\times\partial \Omega$},\\
    u=0, \quad \; \frac{\partial v}{\partial t}=0
    \quad \,  &\mbox{on $\,\{0\}\times \Omega$}\end{array}\right.
 \end{eqnarray}
 \begin{eqnarray}  \left\{ \begin {array}{ll}  -i \frac{\partial v}{\partial
  t}= \sum_{i=1}^n \frac{(1+k|x|^2)^2}{4} \,\,\frac{\partial^2 v}{\partial x_i^2}
  \quad \;\; &\mbox{in $\,(0,+\infty)\times \Omega$}\\
   v=b(t) \quad \,  &\mbox{on $\,(0,
  +\infty)\times\partial \Omega$},\\
    u=0, \quad \; \,  & \mbox{on $\,\{0\}\times \Omega$}, \\
    b(t)\in L^1(0, +\infty) \quad \mbox{and}\;\; \lim_{t\to +\infty}
    b(t)=0.\end{array}\right.
 \end{eqnarray}
  If $v$ satisfies the extra  condition:
   \begin{eqnarray}v(t,x)=a(t), \quad \;
   (t,x)\in (0, +\infty)\times \partial D,\end{eqnarray}
  for some function $a: (0, +\infty)\to (0, +\infty),$
  then $\Omega$ must be a geodesic ball in ${\Bbb M}_k$. }

\vskip 0.29 true cm

\noindent {\it Proof.} \ \  It is easily verified that
  the balance law still holds for the wave equation (5.11) and
  Schr\"{o}dinger's equation (5.12).

  Now, the proof is similar to that of Theorem 1.1, only noticing the following two
techniques: For the wave equation, let us write
\begin{eqnarray} V(s,x)=\sqrt{s} \int_0^{+\infty} v(t, x)e^{-\sqrt{s}\,t}dt,
 \quad \, s>0.\end{eqnarray}
  From (5.11), we get that for any fixed $s>0$,  $V(s,x)$
  satisfies the elliptic boundary value problem:
 \begin{eqnarray} \left\{ \begin{array}{ll} \sum_{i,j=1}^n
 \frac{(1+k|x|^2)^2}{4}\, \frac{\partial^2 V}{\partial
 x_i^2} -s V(s,x) =0 \quad \;\; &\mbox{in}\;\; \Omega,\\
  V=1\quad \;\;  &  \mbox{on}\;\; \partial \Omega. \end{array} \right. \end{eqnarray}
 Moreover, by (5.13) it follows that $V$ is constant on $\partial
 D$. Indeed,
 $$V(s, x)=\sqrt{s}\int_0^{+\infty} a(t)e^{-\sqrt{s}\, t} dt:=c_1(s), \,\quad
 \forall x\in \partial D, \;\; s>0. $$

      For the Schr\"{o}dinger equation (5.12), by putting
 \begin{eqnarray} V(s,x)= \int_0^{+\infty} v(t, x)e^{-ist}dt,
 \quad \, s>0,\end{eqnarray}
 we get elliptic equation
\begin{eqnarray} \left\{ \begin{array}{ll} \sum_{i,j=1}^n
 \frac{(1+k|x|^2)^2}{4}\, \frac{\partial^2 V}{\partial
 x_i^2} -s V(s,x) =0 \quad \;\; &\mbox{in}\;\; \Omega,\\
  V(s,x)=\int_0^{+\infty} b(t)e^{ist} dt\quad \;\;  &
  \mbox{on}\;\; \partial \Omega. \end{array} \right. \end{eqnarray}
 for any fixed $s>0$. Let us replace $V(s,x)$ by $\frac{V(s,t)}{\int_0^\infty b(t)e^{ist}dt}$
 (still denote it by $V(s,x)$), we also obtain the form of (5.15).
 Similarly, we have $$V(s,x)=\int_0^{+\infty} a(t)e^{-ist} dt:=
 c_2(s), \quad \, \forall x\in \partial D,\;\; s>0.$$

 Since $\Omega$ is a bounded domain in ${\Bbb M}_k$ (in the case
  $k\ne 0$, $|x|< \rho$ for any $x\in \bar \Omega$), there exist two
   constants $\alpha>0$ and $\beta>0$ such that
  $$  \alpha\le \frac{1+k|x|^2}{4} \le \beta \quad \;\,
  \text{for all } x\in \bar \Omega.$$
 By using the maximum principle to elliptic equation (5.15), we obtain $V\le 1$ on
 $\bar \Omega$.
Thus, Varadhan's theorem can be applied.
  $\;\;
 \square$

\vskip 1.29 true cm

\bigskip\bigskip
\bigskip \bigskip
\centerline{ACKNOWLEDGEMENT} \bigskip

\thanks{The author would like to thank Professor Rolando Magnanini for the reference \cite{No} and some
 useful suggestions. This work was partially supported by the China Scholarship
 Council (No: 2004307D01) and NSF of China.}


\bigskip




\bigskip\bigskip\bigskip

\begin{thebibliography}{99}





\bibitem{Ah}  L. V. Ahlfors, \textsl{`Old and new in M\"{o}bius groups'},
   Ann. Acad. Sci. Fenn. Ser. A. I. Math. (1984), 93-105.



\bibitem{A}  A. D. Alexandrov,  \textsl{A characteristic property of
spheres}, Ann. Mat. Pura Appl. {\bf 58} (1962), 303-315.



  \bibitem{A1}  G. Alessandrini,  \textsl{Matzoh ball soup:
  a symmetry result for the heat equation}, J. Analyse Math.
   {\bf 54}(1990), 229-236.



  \bibitem{A2}  G. Alessandrini,  \textsl{Characterizing spheres by
  functional relations on solutions of elliptic and parabolic equations},
   Applicable Anal. {\bf 40}(1991), 251-261.



 \bibitem{BBF}    C. Bandle, A. Brillard and M. Flucher, \textsl{Green's function,
 harmonic transplantation, and best Sobolev constant in spaces of
 constant curvature}, Trans. Amer. Math. Soc., {\bf 350}(1998),
  1103-1128.




\bibitem{deB}  N. G. de Bruijn, \textsl{Asymptotic methods in analysis},
 Bibliotheca Mathematica 4, North-Holland Publ. Co., Amsterdam; P.
 Noordhoff Ltd., Groningen; Interscience Publ. Inc., New York, 1958.




\bibitem{Ch1}  I. Chavel, \textsl{Eigenvalues in Riemannian
geometry},  Academic Press, 1984.




\bibitem{Ch2}  I. Chavel, \textsl{Riemannian geometry ---- A
modern introduction}, Cambridge University Press, Second Edition,
2006.



\bibitem{CLN}  B. Chow, P. Lu and L. Ni, \textsl{Hamilton's Ricci flow},
  American Mathematical Society, Providence, Rhode Island, Science Press,
  Beijing, 2006.



\bibitem{DFN}  B. A. Dubrovin, A. T. Fomenko and S. P. Novikov,
 \textsl{Modern geometry ---- Methods and applications, Part I},
 Springer-Verlag, New York, 1984.



\bibitem{Ev}  L. C. Evans, \textsl{Partial differential equations},
Amer. Math. Soc., Providence, RI, 1998.



\bibitem{F1}  A. Friedman, \textsl{Partial differential
equations of parabolic type}, Prentice Hall, Englewood Cliffs, NJ,
1964.



\bibitem{F2}  A. Friedman, \textsl{Partial differential
equations}, Holt, Rinehart and Winston, New York, 1969.


\bibitem{Hi} Noel J. Hicks,  \textsl{Notes on Differential Geometry},
  Van Nostrand, 1965.



\bibitem{Jo}  F. John,
   \textsl{Partial differential equations}, Springer-Verlag, New York Inc,
   Fourth Edition, 1982.



\bibitem{J}  J. Jost, Riemmannian geometry and geometric analysis, Third Edition.
Universitext, Springer-Verlag, Berlin, 2002.



\bibitem{Kl}  M. S. Klamkin, \textsl{ A physical characterization
of a sphere (Problem 64-5*)}, SIAM Review
  {\bf 6}(1964), 61.


\bibitem{Le}   J. M. Lee,  \textsl{Riemannian manifolds}, Springer-Verlag,
New-York Inc., 1997.


\bibitem{Mo1}
  C. B. Morrey, \textsl{On the analyticity of the solutions of analytic
non-linear elliptic systems of partial differential equations},
 Amer. Jour. Math., {\bf 80}(1958), 198-218.



\bibitem{Mo2}  C. B. Morrey,  \textsl{Multiple integrals in the calculus of variations},
 Springer-Verlag, New York, Inc., 1966.




\bibitem{MS1}  R. Magnanini and S. Sakaguchi,
\textsl{The spatial critical points not moving along the heat flow},
    J. Analyse Math. {\bf 71}(1997), 237-261.



\bibitem{MS2}   R. Magnanini and S. Sakaguchi,
\textsl{The spatial critical points not moving along the heat flow
II: The centrosymmetric case}, Math. Z. 230(1999), 695-712,
Corrigendum,  {\it ibid}. {\bf 232}(1999), 389.



\bibitem{MS3}  R. Magnanini and S. Sakaguchi,
\textsl{Matzoh ball soup: Heat conductors with a stationary
isothermic surface}, Ann. of Math. {\bf 156} (2002), 931-946.



\bibitem{No}  J. R. Norris,  \textsl{Heat kernel asymptotics and the distance function
in Lipschitz Riemannian manifolds}, Acta Math. {\bf 179}(1997),
79-103.



\bibitem{Sa1}  S. Sakaguchi, \textsl{When are the spatial level
surfaces of solutions of  diffusion equations invariant with respect
to the time variable}?,  J. Analyse Math. {\bf 78}(1999), 219-243.



\bibitem{Sa2}  S. Sakaguchi, \textsl{Stationary critical points
of the heat flow in spaces of constant curvature},
 J. London Math. Soc. (2) {\bf 63}(2001), 400-412.



\bibitem{Va}  S. R. S. Varadhan, \textsl{On the behavior of
the fundamental solution of the heat equation with variable
coefficients},
  Comm. Pure Appl. Math. {\bf 20}(1967), 430-455.



\bibitem{Wo} J. Wolf,  \textsl{Spaces of constant curvature},
McGraw-Hill, New York, 1967.



\bibitem{Za}  L. Zalcman, \textsl{Some inverse problems of potential theory},
 Contemp. Math. {\bf 63}(1987), 337-350.



\end{thebibliography}
\end{document}